\documentclass{amsart}
\usepackage{amsmath,amsfonts,amssymb}
\usepackage{amsthm}
\usepackage{ifthen}
\usepackage{longtable}
\usepackage{array}
\usepackage{url}
\usepackage{multirow}

\newcommand{\Q}{\mathbb{Q}}

\renewcommand{\P}{\mathbb {P}}

\newtheorem*{theorem*}{Theorem}
\newtheorem{theorem}{Theorem}
\newtheorem{lemma}{Lemma}
\newtheorem{corollary}{Corollary}

\theoremstyle{remark}
\newtheorem{remark}{Remark}

\begin{document}
\title[A note the number of $S$-Diophantine quadruples]{A note on the number of $S$-Diophantine quadruples}
\subjclass[2010]{11D45, 11N32} \keywords{Diophantine equations, $S$-unit equations, $S$-Diophantine tuples}

\author[F. Luca]{Florian Luca}
\address{F. Luca\newline
\indent Mathematical Institute, UNAM Juriquilla\newline
\indent Juriquilla, 76230 Santiago de Quer\'etaro\newline
\indent Quer\'etaro de Arteaga, M\'exico\newline
\indent and\newline
\indent School of Mathematics\newline
\indent University of the Witwatersrand\newline
\indent P. O. Box Wits 2050, South Africa}
\email{fluca@matmor.unam.mx}

\author[V. Ziegler]{Volker Ziegler}
\address{V. Ziegler\newline
\indent Johann Radon Institute for Computational and Applied Mathematics (RICAM)\newline
\indent Austrian Academy of Sciences\newline
\indent Altenbergerstr. 69\newline
\indent A-4040 Linz, Austria}
\email{volker.ziegler\char'100ricam.oeaw.ac.at}

\begin{abstract}
Let $(a_1,\dots, a_m)$ be an $m$-tuple of positive, pairwise distinct integers. If for all $1\leq i< j \leq m$ the prime divisors of $a_ia_j+1$ come from 
the same fixed set $S$, then we call the $m$-tuple $S$-Diophantine. In this note we estimate the number of $S$-Diophantine quadruples in terms of $|S|=r$.
\end{abstract}

\maketitle

\section{Introduction}

There is a vast amount of papers concerning the problem of determining the number of prime divisors of products of the form
$$
\prod_{a \in A, b \in B} (a+b) \qquad \text{and} \quad \prod_{a \in A, b \in B} (ab+1),
$$
where $A$ and $B$ are finite sets of positive integers. In particular, the first product has been studied, first by Erd\H{o}s and Tur\'an \cite{Erdos:1934} 
and their investigations were continued in a series of papers by S\'ark\"ozy and Stewart (see e.g. \cite{Sarkozy:1986,Sarkozy:1994}).
The second product was studied e.g. by Gy\H{o}ry, S\'ark\"ozy and Stewart \cite{Gyory:1996}, S\'ark\"ozy and Stewart \cite{Sarkozy:2000}, and others. 

In their paper \cite{Gyory:1996}, Gy\H{o}ry, S\'ark\"ozy and Stewart conjectured that the largest prime factor of
$$
(ab+1)(ac+1)(bc+1), \qquad 0<a<b<c
$$
goes to infinity as $c$ does. This conjecture has been proved by Corvaja and Zannier~\cite{Corvaja:2003} and Hernandez and Luca 
\cite{Hernandez:2003}, independently. Due to the application of the Subspace theorem their results stay ineffective. The best approach to estimate the growth 
rate of the 
largest prime factor of $(ab+1)(ac+1)(bc+1)$ is due to Luca~\cite{Luca:2005}, who proved that for every fixed finite set of primes $S$, there exist 
ineffective const\-ants $C_{S}$ and $C_{S}'$ such that
$$((bc+1)(ac+1))_{\bar{S}}>\exp\left(C_{S} \frac{\log c}{\log\log c}\right)$$
whenever $a<b<c$ with $c>C_{S}'$, where $(\cdot)_{\bar{S}}$ denotes the $S$-free part.

In case of quadruples effective results are known. For example, Stewart and Tijdeman \cite{Stewart:1997}, proved that the largest prime factor of
$$\prod_{\substack{a,b \in A, \\ a\neq b}} (ab+1)$$
with $|A|\geq 4$, is at least $C \log \log \max A$, where $C$ is an effective computable constant.

Let $S$ be a fixed, finite set of primes. In view of classical Diophantine $m$-tuples we call an $m$-tuple $(a_1,\dots,a_m)$ of positive, pairwise distinct, 
integers $S$-Diophantine if for all $1\leq i <j \leq m$ the set of prime divisors of $a_ia_j+1$ is contained in $S$. The results of Corvaja, Zannier 
\cite{Corvaja:2003} and Hernandez, Luca  \cite{Hernandez:2003} yield the finiteness of $S$-Diophantine triples for fixed $S$. Although we are able to estimate 
the number of $S$-Diophantine triples due to a result of Bugeaud and Luca \cite{Bugeaud:2004a}, it is in principle not possible to determine all triples with 
the methods currently available.

In contrast to the case of triples we can, in principle, effectively determine all $S$-Diophantine quadruples for a given set $S$ due to the result 
of Stewart and Tijdeman \cite{Stewart:1997}. Recently, Szalay and Ziegler \cite{Szalay:2013b}, established an efficient algorithm 
to determine all $S$-Diophantine quadruples for a given set $S$ of primes, provided $|S|=2$. In particular, the results  of Szalay and Ziegler 
\cite{Szalay:2013,Szalay:2013a,Szalay:2013b}, suggest that for $|S|=2$ no quadruple exists at all.

The aim of this note is to give upper bounds for the number of $S$-Diophantine quadruples for fixed sets $S$ of $r$ primes. We need the following 
notations. Let $\Gamma$ be a multiplicative  
subgroup of $\Q^*$ of rank $r$ and denote by $A(n,r)$ an upper bound for the number of non-degenerate solutions $(x_1,\ldots,x_n)\in \Gamma^n$ to the linear 
$S$-unit equation
\begin{equation}\label{eq:UnitEq}
a_1x_1+\dots+a_nx_n =1, \quad a_i \in \Q^*.
\end{equation} 
We call a solution to \eqref{eq:UnitEq} non-degenerate if no subsum on the left hand side of equation \eqref{eq:UnitEq} vanishes. With 
this notation at hand our main result is:

\begin{theorem}\label{Th:quadruples}
Let $S$ be a set of $r$ primes. Then there exist at most 
$$(A(5,r)+A(2,r)^2)A(3,r)$$
$S$-Diophantine quadruples. If $r=2$ or $2\not \in S$, then there exist at most $$A(5,r)A(3,r)$$
$S$-Diophantine quadruples. 
\end{theorem}

Using the best estimates for $A(n,r)$ currently available we obtain

\begin{corollary}\label{Cor:quadruples}
Let $S$ be a set of $r$ primes. Then there exist at most
$$\exp(27398+5126r)$$
$S$-Diophantine quadruples. 
\end{corollary}

In the next section we prove Theorem \ref{Th:quadruples} and in the third section we discuss the number of solutions to  the $S$-unit equation 
\eqref{eq:UnitEq} 
and establish Corollary \ref{Cor:quadruples}.

\section{A system of $S$-unit equations}

Assume that $(a,b,c,d)$ is an $S$-Diophantine quadruple, with $a<b<c<d$. We write,
\begin{align*}
ab+1=&s_1,&ac+1=&s_2, &ad+1=&s_3, \\
bc+1=&s_4,&bd+1=&s_5, &cd+1=&s_6.
\end{align*}
With these notations we have
\begin{align*}
abcd=& s_1s_6-s_1-s_6+1 \\
= & s_2s_5-s_2-s_5+1 \\
= & s_3s_4-s_3-s_4+1
\end{align*}
and obtain the following system of $S$-unit equations
\begin{equation}\label{eq:SUnitsystem}
\begin{split}
s_1s_6-s_1-s_6-s_2s_5+s_2+s_5=&\,\,0,\\
s_1s_6-s_1-s_6-s_3s_4+s_3+s_4=&\,\,0.
\end{split}
\end{equation}
Let us consider the first equation more closely and write $y_1=s_1s_6$, $y_2=s_1$, $y_3=s_6$, $y_4=s_2s_5$, $y_5=s_2$ and $y_6=s_5$. Then the 
first equation of system \eqref{eq:SUnitsystem} takes the form
\[y_1-y_2-y_3-y_4+y_5+y_6=0\]
and has at most $A(5,r)$ projective solutions in $\P^5(\Gamma)$ such that no subsum vanishes, where $\Gamma\subset \Q^*$ is the multiplicative group generated 
by $S$. Note that each projective solution yields at most one solution $(s_1,s_2,s_5,s_6)$.
Indeed, assume $(s_1,s_2,s_5,s_6)$ and $(s'_1,s'_2,s'_5,s'_6)$ correspond to the same projective solution. Then there is a rational number $\rho\neq 0$ such 
that
$s_1=\rho s'_1$, $s_6=\rho s'_6$, $s_2=\rho s'_2$, $s_5=\rho s'_5$ and $s_1s_6=\rho s'_1s'_6$. But this implies that $s_1s_6=\rho^2 s'_1s'_6=\rho s'_1s'_6$, 
thus 
$\rho=1$ and $s_i=s'_i$ for $i=1,2,5,6$.

So we are left to count how many solutions exist with vanishing subsums. Of course there exist no vanishing one-term subsums. Two-term vanishing 
subsums imply either
\begin{itemize}
\item $s_i=s_j$ for $i\neq j$ which is impossible, unless $i,j\in \{3,4\}$ which is excluded, or 
\item $s_i=s_1s_6>abcd>cd+1\geq s_6\geq s_i$ for some $i\in\{1,2,5,6\}$ which is also a contradiction, or
\item $s_i=s_2s_5>abcd>cd+1\geq s_6\geq s_i$ for some $i\in\{1,2,5,6\}$ which is also a contradiction, or
\item $s_1s_6=s_2s_5$, which implies $ab+cd+2=s_1+s_6=s_2+s_5=ac+bd+2$; hence, $(c-b)(d-a)=0$; i.e., $d=a$ or $b=c$, again a contradiction. 
\end{itemize}
Therefore no two-term subsums vanish. Since four- and five-term vanishing subsums imply the existence of two- and one-term vanishing subsums, respectively, we 
are left with the case of three-term vanishing subsums.

Without loss of generality we may assume that the vanishing three-term subsum contains $s_1s_6$. Thus we distinguish whether $s_2s_5$ is contained in the 
vanishing subsum or not. Let us consider the case that $s_2s_5$ is not contained. Then we have an equation of the from $s_1s_6=\pm s_i \pm s_j$. Since 
$s_1=ab+1>2\cdot 1+1>2$ we have $s_1s_6>2s_6>s_i+s_j$ and this case yields no solution.

Therefore both $s_1s_6$ and $s_2s_5$ are contained in the same vanishing three-term subsum and we are left with four systems of $S$-unit equations namely
\begin{equation}\label{eq:3term}
\begin{split}
  s_1s_6-s_5s_2=s_1,& \quad \text{and} \quad s_6=s_5+s_2\\
  s_1s_6-s_5s_2=s_6,& \quad \text{and} \quad s_1=s_5+s_2\\
  s_1s_6-s_5s_2=-s_2,& \quad \text{and} \quad s_1+s_6=s_5\\
  s_1s_6-s_5s_2=-s_5,& \quad \text{and} \quad s_1+s_6=-s_2.
\end{split}
\end{equation}
Note that only the first equation of \eqref{eq:3term} is possible since by assumption $s_1<s_2<s_5<s_6$.
Let $y_1=s_1s_6$, $y_2=s_5s_2$ and $y_3=s_1$. Then the $S$-unit equation 
$$y_1-y_2=y_3$$
has at most $A(2,r)$ projective solutions $(y_1,y_2,y_3)\in\P^2(\Gamma)$. Note that all solutions that yield $S$-Diophantine quadruples are non-degenerate, 
since a vanishing subsum would imply either $s_1s_6=0$ or $s_2s_5=0$ or $s_1=0$.
Each projective solution yields only one possibility for $s_6$. Indeed, assume that $(s_1,s_2,s_5,s_6)$ and $(s'_1,s'_2,s'_5,s'_6)$ yield the same projective 
solution. Then there exists $\rho\in\Q^*$ such that $s_1s_6=\rho s'_1 s'_6=s_1s'_6$, since $s_1=\rho s'_1$, i.e. $s_6=s'_6$. We have now at most $A(2,r)$ 
possible values  for $s_6$; i.e., we are reduced to at most $A(2,r)$ equations of the form 
$$a=s_5+s_2$$
with $a=s_6\neq 0$ fixed. Thus, system \eqref{eq:3term} yields at most $A(2,r)^2$ solutions.

In view of the second statement of Theorem \ref{Th:quadruples} we note that any equation of system \eqref{eq:3term} cannot have a solution if $2\not\in S$. 
Otherwise $s_6$ is odd but $s_5+s_2$ would be even. In case of $r=2$, this implies $S=\{2,p\}$
and the equation  $s_6=s_5+s_2$ turns into 
\begin{equation}\label{eq:Catalanlike}2^{\alpha_6}p^{\beta_6}=2^{\alpha_5}p^{\beta_5}+2^{\alpha_2}p^{\beta_2}.\end{equation}
Considering $2$-adic and $p$-adic valuations, equation \eqref{eq:Catalanlike} reduces to the Diophantine equation
$$ 2^x-p^y=\pm 1.$$
By Mih\u ailescu's solution of Catalan's equation \cite{Mihailescu:2004}, only $p=3$ is possible. On the other hand, Szalay and Ziegler 
\cite{Szalay:2013b} showed that no $\{2,3\}$-Diophantine quadruple exists.

Altogether, we have proved the following result.

\begin{lemma}\label{lem:oneEq}
The first $S$-unit equation in \eqref{eq:SUnitsystem} has at most $A(5,r)+A(2,r)^2$ solutions. If $r=2$ or $2\not \in S$, then there exist at most $A(5,r)$ 
solutions.
\end{lemma}

Now, we turn to the second equation of system \eqref{eq:SUnitsystem}. By Lemma \ref{lem:oneEq}, the first equation in \eqref{eq:SUnitsystem} yields at most 
$A(5,r)+A(2,r)^2$ or $A(5,r)$ many possibilities for the pair $(s_1,s_6)$ respectively. Thus, we assume that the second equation of system 
\eqref{eq:SUnitsystem} is 
of the form
\begin{equation}\label{eq:3termSunit} s_3s_4-s_3-s_4=a \quad \text{with $a\in \Q$ fixed}.\end{equation}
But $S$-unit equation \eqref{eq:3termSunit} has at most $A(3,r)$ solutions provided $a\neq 0$. Indeed no degenerate solution exists since a vanishing subsum 
on the left side of equation \eqref{eq:3termSunit} would imply either
\begin{itemize}
\item $s_3s_4=s_3$ and therefore $s_4=1$, or
\item $s_3s_4=s_4$ and therefore $s_3=1$, or
\item $s_3+s_4=0$ and therefore $s_3s_4<0$.
\end{itemize}

Let us note that $a=s_6s_1-s_6-s_1>2s_6-s_6-s_1>0$, and therefore we have proved the following lemma.

\begin{lemma}
The Diophantine system \eqref{eq:SUnitsystem} has at most $(A(5,r)+A(2,r)^2)A(3,r)$ solutions. If $r=2$ or $2\not \in S$, then there exist at most 
$A(5,r)A(3,r)$ solutions.
\end{lemma}

In order to prove Theorem \ref{Th:quadruples} it remains to prove that for fixed integers $s_1,\dots,s_6$ there exists at most one quadruple $(a,b,c,d)$. Since
\begin{align*}
a&=\sqrt{\frac{(s_1-1)(s_2-1)}{s_4-1}},& b&=\sqrt{\frac{(s_1-1)(s_4-1)}{s_2-1}}, \\
c&=\sqrt{\frac{(s_2-1)(s_4-1)}{s_1-1}},& d&=\sqrt{\frac{(s_5-1)(s_6-1)}{s_4-1}},
\end{align*}
the proof of Theorem \ref{Th:quadruples} is complete.

\section{Proof of Corollary \ref{Cor:quadruples}}

A look through the vast literature on $S$-unit equations shows that for $S$-unit equations over the rationals the best result is due to Evertse 
\cite{Evertse:1984a} provided $|S|=2$ and due to Amoroso and Viada \cite{Amoroso:2009} in the general case. Therefore we may assume $A(2,r)=3\cdot7^{3+2r}$ and 
$A(n,r)=(8n)^{4n^4(n+r+1)}$. A look at the 
proof of the bound for $A(n,r)$ in \cite{Amoroso:2009} shows that this bound is derived by the recursive relation
\[A(n,r)\leq 2^n A(n-1,r)B(n,r+1),\]
where $B(n,r)=(8n)^{6n^3(n+r)}$. Note that this recursive estimate already appears in~\cite{Evertse:2002}. However, recursively 
computing $A(n,r)$ we obtain
$$A(3,r)\leq 8\cdot 3\cdot7^{3+2r} \cdot 24^{162(4+r)}< \exp(2069+518.8 r).$$
Continuing these computations we arrive at
$$A(5,r)< \exp(25329+4616.3r).$$
With these numbers plugged into Theorem \ref{Th:quadruples}, we obtain Corollary \ref{Cor:quadruples}.

\begin{remark}
Let us note that directly applying the bounds due to Evertse \cite{Evertse:1984a} and Amoroso and Viada \cite{Amoroso:2009} would yield the slightly worse 
bound $\exp(73801+15378r)$ for the number of $S$-Diophantine quadruples. A closer inspection of the computation of the quantity $B(n,r)$ due to Amoroso and 
Viada \cite{Amoroso:2009} and Evertse et.al. \cite{Evertse:2002} would further improve the bounds also in view of the new improvements of the Subspace Theorem 
due to Evertse and Feretti \cite{Evertse:2013}. But we are afraid that the gain is too small for such an effort. 
\end{remark}

\section*{Acknowledgement}

We thank the anonymous referee for his valuable suggestions on treating the three-term vanishing subsums in the proof of Theorem \ref{Th:quadruples}.

The first author worked on this paper in Fall of 2013 as a long term guest of  the Special Semester on Applications of Algebra and Number Theory at the 
RICAM, Linz, Austria. He thanks Arne Winterhof for the invitation to participate in this program and RICAM for hospitality. The first author was also supported 
in part by Projects PAPIIT IN104512, CONACyT Mexico--France 193539, CONACyT Mexico--India 163787, and a Marcos Moshinsky Fellowship. The second author was 
supported by the Austrian Science Fund (FWF) under the project P~24801-N26.

\def\cprime{$'$}


\end{document}